\documentclass[11pt,a4paper,authoryear,times]{article}
\usepackage[]{graphicx,float,latexsym,times}
\usepackage{amsfonts,amstext,amsmath,amssymb,amsthm}

 \RequirePackage[OT1]{fontenc}
\long\def\symbolfootnote[#1]#2{\begingroup%
\def\thefootnote{\fnsymbol{footnote}}\footnote[#1]{#2}\endgroup}

\usepackage{amssymb}
\usepackage{xspace}
\usepackage{geometry}
\usepackage{bm}
\usepackage{natbib}
\usepackage{bbm}
\usepackage{color}

\theoremstyle{plain}
\newtheorem{theorem}{Theorem}[section]

\theoremstyle{definition}

\newtheorem{remark}[theorem]{Remark}
 
\newcommand{\beqn}{\begin{eqnarray*}}
\newcommand{\eeqn}{\end{eqnarray*}}
  
\def\ee1{\textrm{\mathversion{bold}$\mathbf{\varepsilon}$\mathversion{normal}}}

\newcommand{\bfz}{{\bf z}}

\def\eu{\mathbf{{u}}}

\newcommand{\N}{\mathbb{N}}
\newcommand{\R}{\mathbb{R}}
\newcommand{\PP}{\mathbb{P}}

\def\XX{\textrm{\mathversion{bold}$\mathbf{X}$\mathversion{normal}}}
\def\eX{\mathbf{X}}
\def\eeX{\mathbb{X}}
\def\ex{\mathbf{x}}
\def\ey{\mathbf{y}}

\newcommand{\Var}{\mathbb{V}\mbox{ar}\,}
\newcommand{\E}{\mathbb{E}\,}

\def\argmin{\mathop{\mathrm{arg\,min}}} 
\usepackage[utf8]{inputenc}    
\DeclareMathAccent{\widehat}{\mathord}{largesymbols}{"62}
\DeclareMathAccent{\widetilde}{\mathord}{largesymbols}{"65}
\def\pth#1{\left(#1\right)}

\def\cro#1{\left[#1\right]}

\def\eeX{\mathbb{X}}

\def\ebo{\textrm{\mathversion{bold}$\mathbf{\beta^0}$\mathversion{normal}}}

\def\eb{\textrm{\mathversion{bold}$\mathbf{\beta}$\mathversion{normal}}}  
  
\def\ed{\textrm{\mathversion{bold}$\mathbf{\delta}$\mathversion{normal}}}

\def\eU{\textrm{\mathversion{bold}$\mathbf{\Upsilon}$\mathversion{normal}}}

\def\XX{\textrm{\mathversion{bold}$\mathbf{X}$\mathversion{normal}}}

\def\ee1{\textrm{\mathversion{bold}$\mathbf{\varepsilon}$\mathversion{normal}}}

\def\eu{\mathbf{{u}}}
\def\e1{1\!\!1}

\def\XX{\textrm{\mathversion{bold}$\mathbf{X}$\mathversion{normal}}}
\def\eX{\mathbf{X}}
\def\eeX{\mathbb{X}}
\def\ex{\mathbf{x}}
\def\ey{\mathbf{y}}

\def\argmin{\mathop{\mathrm{arg\,min}}} 

\def\hh{ \hspace*{0.5cm}}

\begin{document}
\title{Adaptive elastic-net selection  in a  quantile  model with diverging number of variable groups}

\date{}

\author{GABRIELA CIUPERCA  \footnote{{\it email: Gabriela.Ciuperca@univ-lyon1.fr}}}
\date{}
\maketitle

\begin{center}
{\it Universit\'e de Lyon, Universit\'e Lyon 1, 
CNRS, UMR 5208, Institut Camille Jordan, 
Bat.  Braconnier, 43, blvd du 11 novembre 1918, 
F - 69622 Villeurbanne Cedex, France,}\\

\end{center}


\begin{abstract} 
In real applications of the linear   model, the explanatory variables are very often  naturally grouped,  the most common example being the multivariate variance analysis. In the present paper,   a quantile  model with structure group is considered,  the number of groups can diverge with sample size. We introduce and study the adaptive elastic-net group estimator,  for improving the parameter estimation accuracy. This method allows automatic selection, with a probability converging to one, of significant groups and further the non zero  parameter estimators are asymptotically normal. The convergence rate of the adaptive elastic-net group quantile estimator is also obtained, rate  which  depends on the number of groups. In order to put the estimation method into practice, an algorithm based on the subgradient method is proposed and implemented. The Monte Carlo simulations show that the adaptive elastic-net group quantile estimations are more accurate that other existing group estimations in the literature.  Moreover, the numerical study confirms  the theoretical results and the usefulness of the proposed estimation method.
\end{abstract}
 

 
 
\section{Introduction}
\subsection{\textbf{Motivation and state of the art literature review}}
Classically, for a parametric regression model, under assumption that the model errors   are i.i.d. random variables with mean zero and finite variance, the model parameters   are estimated by Least Square (LS) method. If the model errors  don't satisfy the conditions on mean and variance, then the LS method is not appropriate since it can provide inaccurate estimators, i.e. biased and with a large variance. In this case, a very interesting and robust alternative is   the quantile estimation method, introduced by \cite{Koenker-Bassett-78}. For a complete review on properties of the unpenalized quantile estimators, the reader can see book \cite{Koenker-05}, with  developments in \cite{Koenker-Cher-He-Peng.17}. For a linear quantile model with a very large or diverging (with sample size) variable number, it is necessary to automatically detect the  significant variables, without using hypothesis tests. From where  the idea of penalizing the quantile process with  LASSO or adaptive LASSO type penalties, introduced by \cite{Tibshirani.96} and \cite{Zou.06}, respectively, for the LS loss function. The literature being very large on this topic, we can't give here  an exhaustive list, the research on the LASSO method for quantile model being very active in the last decade. \cite{Wu-Liu-09} consider a quantile linear regression, while \cite{Zou-Yuan-08} study a composite quantile regression, the  adaptive LASSO penalty being chosen in both works. Many interesting results on the automatic selection have been obtained when the number of quantile model parameters diverges with the sample size. Examples include, \cite{Gao-Huang.10}, \cite{Belloni-Chernozhukov.11}, \cite{Wang.Wu.Li.12}, \cite{Zheng-Gallagher-Kulasekera-13}, \cite{Zheng-Peng-He.15}. However, quite often in practical applications, grouped explanatory variables must be considered in a linear model, an example being the variance analysis model, with  the goal of the automatic detection of variable groups   that influence the response variable. From where  the  subject of this paper on the detection of the  significant explanatory variable groups in a quantile model when the number of groups depends on the sample size. An adaptive  LASSO  penalty can be considered, to which one will add another penalty to improve the  estimator accuracy. Inspired  by \cite{Zou-Hastie-05}, where  for a linear model with ungrouped explanatory variables was considered, we take the elastic-net as  second penalty.   \\
Previously, for a linear model with the number of ungrouped explanatory variables depending on the number $n$ of samples and  of order $n^c$, with $0 \leq c <1 $,   \cite{Zou-Zhang-09} show the oracle property for  adaptive elastic-net estimator. Recently, always for the LS loss function, \cite{Xin.Hu.Liu.17} generalize this result for multivariate response variable and considering the model with grouped variables. For a model with ungrouped variables, in fixed number, \cite{Slawski.Castell.Titz.10} penalize a convex loss function with an elastic-net penalty, while \cite{Zeng-Xie-14} penalize the LS process with the SCAD function and the $L_2$ norm on the regression coefficients. For a linear model with grouped explanatory variables, the LS loss function is penalized with $L_1$ and $L_2$ penalties in \cite{Simon-Friedman-Hastie-Tibshirani.13}, where  an algorithm is also proposed in order  to fit the model, via accelerated generalized gradient descent and  implement this algorithm in the R package \textit{SGL}.  It should be noted that, algorithms to estimate a quantile model with penalty of elastic-net type, but not with grouped variables, can be found in the following papers: \cite{Slawski.12} develops solution paths algorithms, while, \cite{Yi-Huang-17} propose  an algorithm semismooth Newton coordinate descent,  implemented in the R package \textit{hqreg}.

\subsection{\textbf{Contribution of present paper}}
In the present work, we consider a linear quantile model with grouped explanatory variables, the number of groups  can diverge with the  sample size. Besides to automatically selecting the groups of significant variables, which is usually done using adaptive LASSO penalty, we propose to improve the accuracy of the parameter estimators. For this,  we propose and study the adaptive elastic-net group quantile estimation method. We also show the asymptotic normality for non-zero parameter estimators. For a model with diverging number of variable groups with adaptive elastic-net penalty for a loss function, only the LS loss function has been considered (\cite{Xin.Hu.Liu.17}).  An important contribution of the present paper is to propose an algorithm based on the subgradient method   for computing the adaptive elastic-net group quantile estimator of the regression parameters. Using the  proposed algorithm, simulations are realized to confirm the theoretical results and  to compare our  estimator with existing estimators  for a group quantile linear model.  Point out that, the simulations confirm that adaptive elastic-net group quantile estimator proposed in the present work is more accurate than the adaptive group LASSO quantile estimator existing in literature. Moreover, for a quantile model with ungrouped explanatory variables, the results obtained by simulations using the proposed algorithm are compatible with those obtained by an existing algorithm. \\

The remainder of this paper is organized as follows. Section \ref{Section_Model} introduces the group quantile model,  presents general notations and assumptions. In Section \ref{Section_elastic}, the adaptive elastic-net group quantile estimator is proposed and studied, the convergence rate and oracle property are shown.  An algorithm for computing numerically  the parameter estimations and a criterion for choosing the tuning parameters are proposed in Section \ref{Algorithm}. Section \ref{Simulations} presents numerical results on simulations in order to compare our estimation method with existing adaptive quantile methods.  All proofs are postponed in Section \ref{Section_Proofs}. 
\section{Model}
\label{Section_Model}
Let us consider a linear model with  $g$ groups of explanatory  variables. The number $g$ can depend on $n$, while, the number of explanatory variables in each group don't depend on $n$.  We suppose, without reducing the generality, that each group contains the same number of variables $p$. Let us denote by $r_n = g p$. 
 Thus, the following model is considered in the present paper:
 \begin{equation}
 \label{M1}
 Y_i=\sum^g_{j=1}\eX_{i,j}^t \eb_j+\varepsilon_i=\eeX_i^t \eb^g+\varepsilon_i, \qquad i=1, \cdots, n,
 \end{equation}
 with $\eb^g \equiv (\eb_1, \cdots , \eb_g) \in \R^{r_n}$ and $\eb_j \in \R^p$ the vector of parameters for the group $j$,  for  $j=1, \cdots, g$. For each observation $i$, the vector  $\eeX_i \in \R^{r_n}$ contains all explanatory variables and $\eX_{i,j}$ the explanatory variables of the $j$th group. The explanatory variables  $\eeX_i=(\eX_{i,1}, \cdots , \eX_{i,g})$, with $\eX_{i,j} \in \R^p$,  are assumed to be deterministic, for $j=1, \cdots , g$ and $i=1, \cdots , n$. In model (\ref{M1}), $Y_i$ is the response variable, $\varepsilon_i$ the model error and $\ebo=(\eb^0_1, \cdots , \eb^0_g)$ the true value (unknown) of the parameter $\eb^g$.  For $p=1$,   the model becomes with ungrouped explanatory variables.\\
\hh  For a fixed quantile index $\tau \in (0,1)$, we consider the quantile loss function $\rho_\tau(.): \R \rightarrow(0, \infty)$ defined by $\rho_\tau(u)=u(\tau-\e1_{u <0})$.  
The associated quantile process for model (\ref{M1}) is:
 \begin{equation}
 \label{Gn}
 G_n(\eb^g) \equiv \sum^n_{i=1} \rho_ \tau(Y_i -\eeX_i^t \eb^g).
 \end{equation}
The group quantile estimator is by definition:
 \begin{equation}
 \label{qest}
 \widetilde{{\eb}^{g}}\equiv \argmin_{\eb^g \in \R^{r_n}}G_n(\eb^g). 
 \end{equation}
 For the particular case $\tau=1/2$  we obtain the median regression, for which the quantile process and associated estimator (\ref{qest}), are reduced to the absolute deviation process and to the least absolute deviation estimator, respectively.\\

 For model (\ref{M1}), the purpose of this paper   is to automatically detect  the significant groups of the  explanatory variables, simultaneous with improving the accuracy of the estimators. An adaptive elastic-net method is proposed and studied.  The oracle property is usually requested when the parameter number   of a model diverges with $n$. 
  For a group linear model, the oracle property is the satisfaction of the following two properties: the significant groups of explanatory variables are automatically detected with a probability converging to 1 as $n \rightarrow \infty$ and the non-zero parameters are estimated by asymptotically normal estimators. \\
  Consider  the following set:
 \[  {\cal A}^0 \equiv \{ j \in \{1 , \cdots , g \}; \; \eb^0_j \neq \textbf{0}_p \} ,\]
which  contains the indexes of the significant group of explanatory variables. The set ${\cal A}^0$ is unknown. 
  
\subsection{\textbf{Notation}}
Give now some general notations. All vectors and matrices are denoted by bold symbols and all vectors are column.   For a vector $\textbf{v}$, we denote by $\|\textbf{v}\|_2$  its euclidean norm. We denote also by $\textbf{v}^t$ its transpose and by $|\textbf{v}|$ a vector of the same dimension as $\textbf{v}$, with the components the absolute values of the components of $\textbf{v}$. We use also the notation for the following $p$-dimensional vectors: ${\bf1}_p \equiv(1, \cdots , 1)$ and ${\bf0}_p \equiv(0, \cdots , 0)$. If $\textbf{v}$ and $\textbf{w}$ are two vectors of the same dimension, then, for brevity reason, we denote by $\textbf{v} < \textbf{w}$ the fact that each component of $\textbf{v}$ is less than the corresponding component of $\textbf{w}$. For a positive definite matrix, we use $\mu_{min}(.)$ and $\mu_{max}(.)$ to denote its  smallest and largest eigenvalues. When it is not specified, the convergence is for $n \rightarrow \infty$. Throughout the paper, $C$ denotes a positive generic constant not depend on $n$, which may take a different value in different formula or even in different parts of the same formula. The value of $C$ is not of interest. The notation $\e1_{(.)}$ is used for the indicator function.
On the other hand, let us consider $\varepsilon$ the generic variable for the sequence $(\varepsilon_i)_{1 \leqslant i \leqslant n}$. For an index  set   ${\cal A}$, we denote by $|{\cal A}|$ its cardinality and by $\overline{{\cal A}}$ its complementary set. We also denote by $\eb_{\cal A}$ the subvector of $\eb^g$ containing all subvectors $\eb_j$, with $j \in {\cal A}$.\\
 For two positive sequences $(a_n)$, $(b_n)$ we denote by $a_n \gg b_n$ the fact that $\lim_{n \rightarrow \infty} a_n/b_n=\infty$. Moreover, for two positive sequences  $(a_n)$, $(b_n)$, we write $a_n=O(b_n)$, if $a_n\leq C b_n$ for some constants $C>0$ and for $n$ large enough. We also write $a_n=o(b_n)$ if $a_n/b_n \rightarrow 0$ as $n \rightarrow\infty$. We also use the following notations for two random variable sequences $(U_n)$ and $(V_n)$: $V_n=O_{\PP}(U_n)$ if, for any $\epsilon>0,$ there exists a constant $C>0$ such that $\PP[U_n/V_n>C]<\epsilon$  for $n$ large enough, $V_n=o_{\PP}(U_n)$ when $\lim_{n \rightarrow \infty} \PP[U_n/V_n > \epsilon] =0$ for any $\epsilon >0$.

\subsection{\textbf{Assumptions}}
 For the errors $(\varepsilon_i)_{1 \leqslant i \leqslant n}$ and the design $(\eeX_i)_{1 \leqslant i \leqslant n}$ of model (\ref{M1}), the following general assumptions are considered.   
  \begin{description}
  \item  \textbf{(A1)} Random error  terms $(\varepsilon_i)_{1 \leqslant i \leqslant n}$ are i.i.d. with the  distribution function $F:{\cal G}\rightarrow [0,1]$ and $f$ its density function,  such that $\PP[\varepsilon <0]=\tau$. The set ${\cal G} \subseteq \R$ and $0 \in {\cal G}$.  The density function $f$ is continuously,  strictly positive in a neighbourhood of zero and has a bounded first derivative in the neighbourhood of 0.    
 \item \textbf{(A2)} There exist two constants $0 < m_0 \leq M_0 < \infty$, such that  $ m_0 \leq   \mu_{\min} ( n^{-1} \sum^{n}_{i=1} \eeX_{i} \eeX_{i}^t) \leq    \mu_{\max} ( n^{-1} \sum^{n}_{i=1} \eeX_{i} \eeX_{i}^t) \leq M_0 $.  
\item  \textbf{(A3)}  $\pth{ {g}/{n}}^{1/2} \max_{1 \leqslant i \leqslant n} \| \eeX_i\|_2 \rightarrow 0$, as $n \rightarrow \infty$. 
\item  \textbf{(A4)} The number of groups  is such that $g=O(n^c)$, with $0\leq c < 1$. 
\item \textbf{(A5)} Let be  $h_0\equiv  \min_{1 \leqslant j \leqslant | {\cal A}^0|} \| \eb^0_j\|_2 $. There exists $M>0$ such that for $\alpha >(c-1)/2$  we have $n^{-\alpha} h_0 \geq M$. 
 \end{description}
   Assumption (A1) is standard for a quantile model when the number of parameters can depend on $n$ (see for example \cite{Ciuperca-17a},\cite{Wu-Liu-09}), while Assumption (A2), is standard for linear model (see \cite{Wu-Liu-09}, \cite{Zhang-Xiang-15}). 
       Assumptions (A3) and (A4) are considered in models with number of parameters depending on $n$, see \cite{Zhang-Xiang-15}, \cite{Zou-Zhang-09}, \cite{Ciuperca-17a}. 
By assumption (A5), considered also by \cite{Zhang-Xiang-15}, we suppose that for significant groups of explanatory variables, the true values of the parameters can depend on $n$ and converge to 0 as $n \rightarrow \infty$, when $\alpha >0$.  
 
\section{Asymptotic properties}
 \label{Section_elastic}
In this section we consider that the number $g$ of groups is of order $n^c$, with   the constant $c$ such that  $0\leq c<1$. If $c=0$ then the group number is fixed.  It is interesting to note that the   results stated in this  section are original also for a model with ungrouped explanatory variables ($p=1$). \\
Since the model is well defined now, we can give more precise definition of the oracle property for a parameter estimator. For a model with grouped explanatory variables, an estimator satisfy the oracle property if the \textit{estimation sparsity (consistent selection)} and the \textit{asymptotic normality} properties are satisfied simultaneously. More precisely, the  sparsity property (also say, consistency in selection) is when the true non-zero vectors are estimated as non-zero and the null parameter vectors are shrunk directly as a null vector, with a probability converging to 1 as $n \rightarrow \infty$. The asymptotic normality  property is when the estimators of the true non-null    parameter vectors have a normal asymptotic distribution. \\
In order to have the oracle property for the estimators but also to   improve the accuracy of parameter estimators, we propose a new estimator, inspired  by the work of \cite{Zou-Zhang-09}, by penalizing the quantile process with an adaptive weighted penalty and with an elastic-net penalty. 
In \cite{Zou-Zhang-09}, the adaptive elastic-net estimator is considered for $p=1$, $g=n^c$, with $0 \leq c <1$, for a linear regression with LS loss function.   \\ 
\hh For model (\ref{M1}), for $\eb^g \in \R^{r_n}$, with $r_n=gp$, let us define the adaptive elastic-net  group quantile process   by:
\begin{equation}
\label{aeQn}
E_n(\eb^g) \equiv G_n(\eb^g) +\lambda_{1n}\sum^g_{j=1} \widetilde{\omega}_{n;j} \| \eb_j\|_{2} +\lambda_{2n} \sum^g_{j=1} \| \eb_j\|_2^2,
\end{equation}
with the process $G_n(\eb^g)$ defined by (\ref{Gn}), the weights $\widetilde{\omega}_{n;j} \equiv \| \widetilde \eb_j \|_2^{- \gamma}$, the vector $\widetilde \eb_j$  being the $j$-th component of the quantile estimator vector defined by relation (\ref{qest}) and $\gamma>0$ a positive constant which will be later specified. 
For the particular case $\lambda_{2n} =0$, for any $n  \in \N$, we obtain the adaptive quantile process, considered in \cite{Ciuperca-17a}.\\
The adaptive elastic-net group quantile estimator is defined by:
\begin{equation}
\label{ael-net}
\widehat{{\eb}^g}=(\widehat{\eb}_1, \cdots , \widehat{\eb}_g) = \argmin_{\eb^g \in \R^{r_n} } E_n(\eb^g). 
\end{equation}
In a similar way to the set ${\cal A}^0$,  let us consider the following sets of indexes:   
\[
\widehat{{\cal A}}_n  \equiv \{ j \in \{1 , \cdots , g \}; \; \widehat{\eb}_j \neq \textbf{0}_p \}.
\]
The set $\widehat{{\cal A}}_n$ contains the indexes of the non-zero vector estimators. 
In the present work,  the set ${\cal A}^0$     is estimated through the  adaptive elastic-net method by the set  $\widehat{{\cal A}}_n $.\\
 
For studying the asymptotic behaviour $\widehat{{\eb}^g}$ we must first know the asymptotic properties of the group quantile estimator $\widetilde{{\eb}^{g}}$  which intervenes in the weights $\widetilde{\omega}_{n;j}$. More precisely, we need to know if $\widetilde{{\eb}^{g}}$ is consistent and if so, what is its   convergence rate.   
Under assumptions (A1)-(A4), we have by Lemma 1 of \cite{Ciuperca-17a} that the convergence rate of the group quantile estimator $\widetilde {\eb^g}$ is of order $(g/n)^{1/2}$:  $\| \widetilde {\eb^g }-\ebo \|_2=O_{\PP}((g/n)^{1/2}) $.\\
A first studied property for the  adaptive elastic-net group quantile estimator $\widehat{{\eb}^g}$ is the convergence rate. In order to find the convergence rate of the $\widehat{{\eb}^g}$ in the case of $g$ depending on $n$, we suppose for the tuning parameters $\lambda_{1n}$, $\lambda_{2n}$, as $n \rightarrow \infty$, that,
\begin{align}
\label{C1bisb}
&\lambda_{1n} \rightarrow \infty, \quad \lambda_{1n} n^{(c-1)/2- \alpha \gamma} \rightarrow 0.\\
\label{C1bis}
&\lambda_{2n} \rightarrow \infty,  \quad \lambda_{2n}n^{(c-1)/2} \rightarrow 0.
\end{align}
Condition  (\ref{C1bisb}) was also considered by \cite{Ciuperca-17a}  for studying the asymptotic behaviour of  the adaptive group LASSO quantile estimator.\\
By the following theorem we show that the convergence rate of  $\widehat{{\eb}^g}$  is of order $\big(g/n \big)^{1/2}$, which for $c=0$ ($g$ fixed) becomes $n^{-1/2}$, the convergence rate obtained also in the proof of Remark \ref{theorem1_SP}, under more appropriate assumptions. 
\begin{theorem}
\label{th_vconv_SP}
If $0\leq  c< 1$, under assumptions (A1)-(A5) and conditions (\ref{C1bisb}), (\ref{C1bis}) for the tuning parameters, we have   $\| \widehat{{\eb}^g}- \ebo\|_2=O_{\PP}\big((g/n)^{1/2} \big)$.
\end{theorem}

With this result we can now study the oracle property of $\widehat{{\eb}^g}$. For this, we consider the following additional conditions, as $n \rightarrow \infty$, for the tuning parameters: 
\begin{align}
  \label{C1bis2}
 &\displaystyle{\lambda_{1n} n^{(1-c)(1+\gamma)/2-1} \rightarrow \infty }
 \\
 \label{C2bis1}
& \displaystyle{ \lambda_{1n} |{\cal A}^0|  n^{-(1+c)/2 } \rightarrow 0. }
 \\
 \label{C2bis2}
&\displaystyle{  \lambda_{2n} |{\cal A}^0|  n^{-(1+c)/2 } \rightarrow 0. }
 \end{align}
If in assumption (A5) we consider    $\alpha=0$,  that is, the smallest norm of significant groups don't depend on $n$ and if $|{\cal A}^0| < \infty$, then condition (\ref{C1bisb}) implies condition (\ref{C2bis1}). If $|{\cal A}^0| < \infty$, condition (\ref{C1bis}) implies condition (\ref{C2bis2}).\\
 We denote by $\eeX_{i,{\cal A}^0}$ the columns  $\XX_{i,j} $ with $j \in {\cal A}^0$. The following theorem states that the adaptive elastic-net group quantile estimator satisfies the oracle property. 
\begin{theorem}
\label{Theorem 2SPL} If $0\leq  c< 1$, suppose  that assumptions (A1)-(A5)  are  satisfied  and also that the tuning parameters satisfy  (\ref{C1bisb}), (\ref{C1bis}), (\ref{C1bis2}), (\ref{C2bis1}), (\ref{C2bis2}).  Then:\\
(i) $\PP \cro{\widehat{\cal A}_n ={\cal A}^0}\rightarrow 1$,  for $n \rightarrow \infty $.\\
(ii) For any  vector $\eu$ of size $(p |{\cal A}^0|)$ such that $\| \eu\|_2=1$, with notation $\eU_{n,{\cal A}^0} \equiv n^{-1} \sum^n_{i=1} \eeX_{i,{\cal A}^0}  \eeX_{i,{\cal A}^0}^t$, we have,  $n^{1/2} (\eu^t \eU^{-1}_{n,{\cal A}^0} \eu)^{-1/2} \eu^t ( \widehat{{\eb}^g} - \eb^0)_{{\cal{A}}^0}  \overset{\cal L} {\underset{n \rightarrow \infty}{\longrightarrow}} {\cal N}\big(0, \tau (1- \tau ) f^{-2}(0) \big)$.
\end{theorem}
More precisely, conditions (\ref{C1bisb}), (\ref{C1bis}) and (\ref{C1bis2}) are required in the previous theorem for the tuning parameters $\lambda_{1n}$, $\lambda_{2n}$,  for proving the estimator sparsity in Theorem \ref{Theorem 2SPL}\textit{(i)}. Conditions (\ref{C2bis1}), (\ref{C2bis2})  are necessary for the asymptotic normality (see the proof of Theorem \ref{Theorem 2SPL} in Section \ref{Section_Proofs}).  
The asymptotic normality of   Theorem \ref{Theorem 2SPL}\textit{(ii)}  shows that the elastic-net penalty  don't affect the asymptotic law of the non-zero parameter estimators. The Gaussian limit distribution being the same as for a quantile estimator, without penalty (see \cite{Koenker-05}. The oracle property obtained in Theorem \ref{Theorem 2SPL}  has been shown for other models or other penalties when the number of parameters diverges with the sample size. For a quantile linear model with ungrouped  explanatory variables, \cite{Zheng-Gallagher-Kulasekera-13},  \cite{Zheng-Peng-He.15} and \cite{Ciuperca-17a} considered an adaptive LASSO penalty, while \cite{Ciuperca-15} used a seamless $L_0$ penalty. \cite{Wang.Wu.Li.12} shows the sparsity of the quantile estimators with a SCAD and MCP penalty. To the knowledge of the author, the adaptive elastic-net penalty  was considered only for a LS loss function by \cite{Zou-Zhang-09} which shows the oracle property for a model with ungrouped variables and by \cite{Xin.Hu.Liu.17} which shows the starsity property for a model with grouped variables.   \\

Let's take a closer look at the  case $c=0$, in this case the conditions on the design and on the tuning parameters can be simplified. Let us  denote $r=r_n$, since it don't depend on $n$. 
 For the case $g$ fixed, assumption (A3) implies $n^{-1} \max_{1 \leq i \leq n} \eeX_i^t \eeX_i {\underset{n \rightarrow \infty}{\longrightarrow}}  0$.  Then, instead of assumptions (A2), (A3) we consider:
 \begin{description}
  \item  \textbf{(A6)} $n^{-1} \max_{1 \leq i \leq n} \eeX_i^t \eeX_i {\underset{n \rightarrow \infty}{\longrightarrow}}  0$ and $n^{-1} \sum^{n}_{i=1} \eeX_i \eeX_i^t {\underset{n \rightarrow \infty}{\longrightarrow}} \eU$, with $\eU$ a $r \times r$ positive definite matrix.
  \end{description}
 
For the penalties, the tuning parameters $\lambda_{1n}$, $\lambda_{2n}$ and the power $\gamma$ of weight $\widetilde{\omega}_{n;j}$ are such that, for $n \rightarrow\infty$:
\begin{equation}
\label{C1}
\lambda_{1n} \rightarrow \infty,\quad \lambda_{2n} \rightarrow \infty, \quad n^{-1/2} \lambda_{1n} \rightarrow 0,  \quad  n^{(\gamma-1)/2} \lambda_{1n} \rightarrow \infty, \quad  n^{-1/2} \lambda_{2n} \rightarrow 0.
\end{equation}
For the tuning parameter $\lambda_{1n}$ the considered conditions are the same as for the adaptive group LASSO quantile estimator of \cite{Ciuperca-17a}. By assumption $n^{-1/2} \lambda_{1n} \rightarrow 0$ of relation (\ref{C1}) and $n^{\gamma/2-1} \lambda_{1n} \rightarrow \infty$, as $n \rightarrow \infty$,  considered in Remark \ref{th_selection_SP} we deduce that $\gamma >1$.\\

The following remark shows that for $g$  fixed, only assumptions (A1), (A6) and condition (\ref{C1}) for the tuning parameters are needed for the asymptotic normality of $\widehat{\eb^g}_{{\cal A}^0}$. The proof of Remark \ref{theorem1_SP} is given in Section \ref{Section_Proofs}.
\begin{remark}
\label{theorem1_SP}
\textit{If $c=0$, under assumptions (A1), (A6) and condition (\ref{C1}),  we have $n^{1/2}(\widehat{{\eb}^g} - \eb^0)_{{\cal{A}}^0}  \overset{\cal L} {\underset{n \rightarrow \infty}{\longrightarrow}} {\cal N}\big(\textbf{0}_{p|{\cal A}^0|}, \tau(1-\tau) f^{-2}(0)\eU^{-1}_{{\cal A}^0} \big)$, with $(p|{\cal A}^0|)$-squared matrix  $\eU_{{\cal A}^0}$ defined as the submatrix of $\eU$  such that with the indexes of rows and columns in ${\cal A}^0$.}
\end{remark}
We also note that Remark \ref{theorem1_SP} is a particular case of Theorem \ref{Theorem 2SPL}\textit{(ii)} when $|{\cal A}^0|$ is bounded.\\ 
 The following remark gives that the elements of ${\cal A}^0$ and $\widehat{\cal A}_n$ coincide with a probability converging to one as $n \rightarrow \infty$.  
 \begin{remark}
\label{th_selection_SP}
\textit{If $c=0$, under assumptions (A1), (A6) and condition (\ref{C1}) together   $n^{\gamma/2-1} \lambda_{1n} \rightarrow \infty$, as $n \rightarrow \infty$, we have $\lim_{n \rightarrow \infty }\PP [\widehat{\cal A}_n ={\cal A}^0]=  1 $.}
\end{remark}
The proof of Remark \ref{th_selection_SP}, based on the Karush-Kuhn-Tucker  optimality conditions (given in subsection \ref{Algorithm}), is omitted since it is  similar to that of  the proof of Theorem 2 in \cite{Ciuperca-17a}.\\
Always for a quantile model with grouped explanatory variables, but in finite number,  the same   variance matrix of the centered   normal limit distribution  in Theorem \ref{th_selection_SP} was obtained by \cite{Ciuperca-17a} for adaptive LASSO penalty and by \cite{Ciuperca-17b} for adaptive fused LASSO penalty.  The oracle property is also true for  the adaptive LASSO estimator when the  loss function is the LS (see \cite{Wang.Leng.08}. 
 \section{Algorithm}
 \label{Algorithm}
In this section we propose an algorithm, based on the subgradient method, in order to compute the adaptive elastic-net estimator for a group quantile linear model. For this, we will write  the Karush-Kuhn-Tucker (KKT) optimality conditions.\\
 For all  $ j \in \widehat{{\cal A}}_n $,   the following   $p$ equalities  hold with probability one,
\begin{equation}
\label{KKTi}
\tau \sum^n_{i=1} \XX_{i,j} - \sum^n_{i=1} \XX_{i,j} \e1_{Y_i < \eeX^t_i \widehat{{\eb}^g}} - 2 \lambda_{2n} \widehat{\eb}_{j} = \frac{\lambda_{1n} \widetilde{\omega}_{n;j} \widehat{\eb}_{j}}{\| \widehat{\eb}_{j}\|_2}.
\end{equation}
Moreover, for all  $ j \not \in \widehat{{\cal A}}_n  $, for all $k =1, \cdots , p$, we have, with probability one, the following inequality
\begin{equation}
\label{KKTii}
\left| \tau \sum^n_{i=1} X_{i,j,k} - \sum^n_{i=1} X_{i,j,k} \e1_{Y_i < \eeX^t_i \widehat{{\eb}^g}} - 2 \lambda_{2n} \widehat{\beta}_{j,k} \right| \leq \lambda_{1n} \widetilde{\omega}_{n;j} .
\end{equation}
 The $p$ inequalities of relation (\ref{KKTii}), for all  $ j \not \in \widehat{{\cal A}}_n  $, can be also written:
 \[
\bigg| \tau \sum^n_{i=1} \XX_{i,j} - \sum^n_{i=1} \XX_{i,j} \e1_{Y_i < \eeX^t_i \widehat{{\eb}^g}} - 2 \lambda_{2n} \widehat{\eb}_{j}\bigg| \leq   \lambda_{1n} \widetilde{\omega}_{n;j}  |\textbf{s}_j |,
 \] 
 with $\textbf{s}_j$ a $p$-dimensional vector with each component in absolute value less than 1.\\
 From this last relation, together relation (\ref{KKTi}), we consider  the following $g  p$ gradient equations, for any $j \in \{2, \cdots , g \}$: 
 \begin{equation}
 \label{KKT}
 \tau \sum^n_{i=1} \XX_{i,j} - \sum^n_{i=1} \XX_{i,j} \e1_{Y_i < \eeX^t_i  {\eb}^g } - 2 \lambda_{2n}  \eb_{j} +\lambda_{1n} \widetilde{\omega}_{n;j} \textbf{s}_j =\textbf{0}_p,
 \end{equation}
 with $\textbf{s}_j$ a $p$-dimensional vector   such that, 
 \begin{itemize}
 \item if $| \eb_j\|_2 \neq 0$, then $\textbf{s}_j=- \frac{\eb_j}{\| \eb_j\|_2}$, 
 \item if $\|\eb_j\|_2 = 0$, then $ |\textbf{s}_j| \leq {\bf1}_p$. \\
 \end{itemize}
\hh Let us denote by $\eeX_{i,-j}$ the vector $\eeX_i$ without subvector  $\XX_{i,j}$ and $\eb^g_{-j}$ the vector $\eb^g$ without the subvector $\eb_j$. Consequently, if $\|\eb_j\|_2 = 0 $ then relation (\ref{KKT}) becomes:
 $
 \sum^n_{i=1} \XX_{i,j}  (\tau - \e1_{Y_i < \eeX^t_{i,-j}  {\eb}^g_{-j} }) +\lambda_{1n} \widetilde{\omega}_{n;j} \textbf{s}_j =\textbf{0}_p$,
 which implies 
 \[
 \textbf{s}_j =\frac{1}{\lambda_{1n} \widetilde{\omega}_{n;j}}\sum^n_{i=1} \XX_{i,j}  ( \e1_{Y_i < \eeX^t_{i,-j}  {\eb}^g_{-j}  }-\tau).
 \]
Thereby,  if $|\textbf{s}_j| <{\bf1}_p$, then $\big|\sum^n_{i=1} \XX_{i,j}  ( \e1_{Y_i <\eeX^t_{i,-j}  {\eb}^g_{-j} }-\tau)  \big| < \lambda_{1n} \widetilde{\omega}_{n;j}{\bf1}_p$. Thus, in the  algorithm, we consider when $|\textbf{s}_j| <{\bf1}_p$, that is  $\big|\sum^n_{i=1} \XX_{i,j}  ( \e1_{Y_i < \eeX^t_{i,-j}  {\eb}^g_{-j} }-\tau)  \big| < \lambda_{1n} \widetilde{\omega}_{n;j}{\bf1}_p$, that $\eb_j=\textbf{0}_p$. Otherwise:
 \[
 \sum^n_{i=1} \XX_{i,j}  (\tau - \e1_{Y_i < \eeX^t_{i,-j}  {\eb}^g_{-j} })- 2 \lambda_{2n}  \eb_{j}  = \frac{\lambda_{1n} \widetilde{\omega}_{n;j} \eb_{j}}{\| \eb_{j}\|_2},
 \]
 from where, 
 \begin{equation}
 \label{bb}
 \eb_j=\frac{\sum^n_{i=1} \XX_{i,j}  (\tau - \e1_{Y_i < \eeX^t_{i,-j}  {\eb}^g_{-j} })}{2 \lambda_{2n} +\lambda_{1n} \widetilde{\omega}_{n;j}\| \eb_{j}\|_2^{-1}}
 \end{equation}
 and considering the euclidean  norm we obtain: $\|\sum^n_{i=1} \XX_{i,j}  (\tau - \e1_{Y_i < \eeX^t_{i,-j}  {\eb}^g_{-j} }) \|_2=2 \lambda_{2n} \| \eb_{j}\|_2+ \lambda_{1n} \widetilde{\omega}_{n;j}$. From this last relation we get:
 \[
 \| \eb_{j}\|_2 =\frac{\|\sum^n_{i=1} \XX_{i,j}  (\tau - \e1_{Y_i < \eeX^t_{i,-j}  {\eb}^g_{-j} }) \|_2-  \lambda_{1n} \widetilde{\omega}_{n;j}}{2 \lambda_{2n}}
 \] 
and replacing in  relation (\ref{bb}), we obtain:
\[
\eb_j=\frac{\sum^n_{i=1} \XX_{i,j}  (\tau - \e1_{Y_i < \eeX^t_{i,-j}  {\eb}^g_{-j} })}{2 \lambda_{2n} +2\lambda_{1n}\lambda_{2n}\widetilde{\omega}_{n;j}  \big(\|\sum^n_{i=1} \XX_{i,j}  (\tau - \e1_{Y_i < \eeX^t_{i,-j}  {\eb}^g_{-j} }) \|_2-  \lambda_{1n} \widetilde{\omega}_{n;j}\big)^{-1}}.
\] 
Then, from these relations, we can  propose the following algorithm, for $\lambda_{1n}$, $\lambda_{2n}$, $\tau$, $\gamma$ fixed.  \\
\newpage
\textbf{Algorithm}
 \noindent {\bf \hrule }
 \medskip 
\begin{description}
  \item \textit{Step $0$}
\begin{description}
\item set the initial values of $\eb^{(0)}$ for the coefficient parameters. 
\end{description}
\item \textit{Step $k$}\\
For all $j=1, \cdots , g$ we calculate:
\begin{itemize}
\item if $\big|\sum^n_{i=1} \XX_{i,j}  ( \e1_{Y_i < \eeX^t_{i,-j}  {\eb}^{g(k-1)}_{-j} }-\tau)  \big| < \lambda_{1n} \widetilde{\omega}_{n;j}{\bf1}_p$, then, $\eb_j^{(k)}=\textbf{0}_p$,
\item otherwise:
\[
\eb_j^{(k)}=\frac{\sum^n_{i=1} \XX_{i,j}  (\tau - \e1_{Y_i < \eeX^t_{i,-j}  {\eb}^{g(k-1)}_{-j} })}{2 \lambda_{2n} +2 \lambda_{1n} \lambda_{2n}\widetilde{\omega}_{n;j}  \big(\|\sum^n_{i=1} \XX_{i,j}  (\tau - \e1_{Y_i < \eeX^t_{i,-j}  {\eb}^{g(k-1)}_{-j} }) \|_2-  \lambda_{1n} \widetilde{\omega}_{n;j}\big)^{-1}}.
\] 
\end{itemize}
\item \textit{Stopping the algorithm:} The algorithm stops when $\|\eb_j^{(k)} - \eb_j^{(k-1)} \| < \epsilon $, with $\epsilon$ a specified precision.
\end{description}
 \noindent {\bf \hrule }
 \medskip 
\vspace{0.25cm} 
As a starting point, to Step $0$, we can take either the group quantile estimator  $\eb^{(0)} =\widetilde{{\eb}^{g}}$, given by (\ref{qest}), or the adaptive group LASSO quantile estimator proposed by \cite{Ciuperca-17a}. For simulations, we will consider this last estimator.\\

The reader  can find in \cite{Tseng.01} the convergence proprties of a block coordinate descendent methods applied to minimize a nondifferentiable continuous function.\\

 For choosing the tuning parameters $\lambda_{1n}$ and $\lambda_{2n}$, we can use a criterion of type BIC. For this, let us consider $(S_n)_{n \in \N}$ a deterministic sequence defined by:
\begin{itemize}
\item if $g$ is fixed or $g=O(n^c)$ such that $g(\log n)^{-1}=o(1)$, then $S_n = 1$, for any $n \in \N$;
\item if $g=O(n^c)$ such that $g(\log n)^{-1}\neq o(1)$, we take $(S_n)$ converging to $+ \infty$ such that $ g^{-1}(\log n) S_n $ $ {\underset{n \rightarrow \infty}{\longrightarrow}} \infty$ and $ n^{-1}(\log n) |{\cal A}^0| S_n {\underset{n \rightarrow \infty}{\longrightarrow}}  0$.
\end{itemize}
Then, we consider the following criterion :
\begin{equation}
\label{crit}
BIC(\lambda_{1n},\lambda_{2n}) \equiv \log (n^{-1} G_n(\widehat{{\eb}^g}))+\frac{\log n}{n}S_n |\widehat{\cal A}_n|
\end{equation}
and we choose $\lambda_{1n},\lambda_{2n}$ that minimize the criterion.\\
 
Criterion (\ref{crit}) is of type BIC, introduced  by \cite{Wang.Wu.Li.12} for a linear model with diverging number of parameters and estimated by penalized LS method.  For choosing the tuning parameters, the same type of criterion have been proposed for quantile models with ungrouped variables by \cite{Ciuperca-15} for seamless $L_0$ penalty, by \cite{Zheng-Peng-He.15} for adaptive LASSO penalty.  Taking into account the  sparsity  of the parameter  estimator $\widehat{{\eb}^g}$ proved in Theorem \ref{Theorem 2SPL}, the choice of sequence $(S_n)$, is such that it doesn't  allow an overfitted or underfitted model. The condition $ n^{-1}(\log n) |{\cal A}^0| S_n {\underset{n \rightarrow \infty}{\longrightarrow}} 0$ avoids an overfitted model and  $(\log n) S_n g^{-1}  {\underset{n \rightarrow \infty}{\longrightarrow}} \infty$ avoids an underfitted model.
\begin{remark}
\textit{In practice, the choice of quantile index can be done as follows. For the explained variable we calculate  the standardized values: $\widetilde{y}_i=(y_i- \overline{y}_n)/\widehat{\sigma}_y$, with $\overline{y}_n$, $\widehat{\sigma}_y$ empirical mean and standard deviation of $Y$, respectively. Afterwards, we calculate the empirical estimation of $\tau$ by $\widehat{\tau}_n=n^{-1} \sum^n_{i=1}\e1_{\widetilde{y}_i <0}$.}
\end{remark}
 
\section{Simulation study}
\label{Simulations}
In this section we conduct Monte Carlo simulations in order to evaluate and compare our proposed adaptive elastic-net estimator with existing adaptive quantile estimators. The simulations are based on the algorithm presented in subsection \ref{Algorithm}.\\
The considered design is such that:
\[
X_{p(j-1)+k}=\frac{Z_j+R_{p(j-1)+k}}{\sqrt{2}}, \qquad 1 \leq j \leq g, \quad 1 \leq k \leq p,
\]  
with $Z_j$  multivariate normal distribution of mean zero and covariance $Cov(Z_{j_1}, Z_{j_2})=0.6^{| j_1-j_2|}$. Moreover, $R_1, \cdots , R_{r_n}$ are independent standard normal variables. The considered quantile index is $\tau=0.5$. In fact, the design is similar to the one considered by \cite{Ciuperca-17a} for adaptive group LASSO quantile method and by \cite{Wei-Huang-10} for adaptive group LASSO-LS model.  In all simulations, the number of significant groups of  explanatory variables was assumed to be four: $|{\cal A}^0|=4$. The number of  non-significant groups of explanatory variables will be varied. \\
In each simulation, for a data set of dimension $n$, the design and the error distributions are generated. Each simulation is repeated 1000 times. 
\subsection{Ungrouped variables}
In this subsection we consider ungrouped variables. We will  compare our estimation method calculated by algorithm in subsection \ref{Algorithm}, with those obtained by R package \textit{hqreg} for a quantile model with adaptive elastic-net penalty (of \cite{Yi-Huang-17}) and with those obtained by R package \textit{quantreg} for quantile model with adaptive LASSO penalty.\\
We consider four significant variables, with the true values of the parameters $\beta^1_0=0.5$, $\beta^2_0=1$, $\beta^3_0=-1$, $\beta^4_0=-1.5$.
For the  adaptive weights in the elastic-net and LASSO penalties, we considered the power $\gamma=12.25/10$. For the adaptive LASSO method, the considered tuning parameter is $n^{2/5}$ and for adaptive elastic-net estimator  in the R package \textit{hqreg}, the considered tuning parameter is $n^{9/20}$. For the adaptive elastic-net quantile method of the present paper we consider the following tuning parameters: $\lambda_{1n}=n^{1-\gamma/2+1/n}$ and $\lambda_{2n}=c_1n^{2/5}$, with $c_1$ varied on a value grid such that criterion (\ref{crit}) becomes minimal. For two gaussian errors, the results are given Table \ref{Tabl1} where are presented: the median of $|\widehat{\cal A}_n|$, that is the median of  the estimated number of groups with significant explanatory variables by three adapted penalized method estimation: adaptive LASSO quantile (\textit{ag}), adaptive elastic-net quantile (\textit{aEq}) of \cite{Yi-Huang-17} and adaptive elastic-net group quantile (\textit{aEGq}) method  proposed in the present paper, subsection \ref{Algorithm}. We also give in Table \ref{Tabl1} the median of $|\overline{\widehat{\cal A}_n}|$ of the estimated number of non-significant groups by the three estimation methods. Always for the three methods, we calculate the standard-deviation of the parameter estimations and the mean of the absolute value of $Y-\hat Y$, with $\hat Y$ the corresponding model prediction of $Y$.  From the significant and non significant group identification point of view, the results are similar by the three estimation methods. The \textit{aq} and \textit{aEq} methods   provide parameter estimations   with the same standard-deviation, while by the \textit{aEGq} method, the parameter estimations are more precise. For the prevision of  the response variable $Y$, for small $n$, the prevision by \textit{aEGq} method is slightly worse, but when $n$ increases, we obtain the same precisions by \textit{aEGq} and \textit{aEg} methods.

\begin{table} 
{\scriptsize
\caption{\footnotesize  Simulation results for models with ungrouped variables: $p=1$. Estimation methods: adaptive   LASSO quantile \textit{(aq)}, adaptive elastic-net quantile \textit{(aEq)},  adaptive elastic-net group quantile \textit{(aEGq)}. }
\begin{center}
\begin{tabular}{|c|c|c|c|c|c|c|c|c|c|c|c|c|c|c|}\hline
 & & &  \multicolumn{3}{c|}{median true $\neq 0$} & \multicolumn{3}{c|}{median true $= 0$} &  \multicolumn{3}{c|}{sd $(\beta^0-\hat \beta)$} & \multicolumn{3}{c|}{$mean|Y-\hat Y|$}\\
  \cline{4-15} 
 n &g & $\varepsilon$ & aq & aEq & aEGq & aq & aEq& aEGq & aq & aEq & aaEGq & aq & aEq & aEGq  \\ \hline
 30 & 5 & ${\cal N}(0,3)$ & 2 & 3 & 2 & 1 & 0 & 1 & 1.33 & 1.47 & 0.96 & 2.25 & 2.18 & 2.61 \\ 
  & & ${\cal N}(0,1)$ & 3 & 4 & 3 & 1 & 1 & 1 & 1.15 & 1.27 & 0.97 & 0.86 & 0.75 & 1.15 \\ \hline
  60 & 5 & ${\cal N}(0,3)$ & 3 & 3 & 2 & 1 & 1 & 1 & 1.28 & 1.30 & 1.20 & 2.33 & 2.30 & 2.50 \\
  & & ${\cal N}(0,1)$ & 4 & 4 & 3 & 1 & 1 & 1 & 1.26 & 1.28 & 1.29 & 0.79 & 0.78 & 0.90 \\ \hline
  100 & 5 & ${\cal N}(0,3)$ & 3 & 3 & 3 & 1 & 1 & 1 & 1.24 & 1.22 & 1.07 & 2.37 & 2.39 & 2.55 \\
  & & ${\cal N}(0,1)$ & 4 & 4 & 4 & 1 & 1 & 1 & 1.26 & 1.25 & 1.29 & 0.80 & 0.82 & 0.84 \\ \hline
  100 & 20 & ${\cal N}(0,3)$ & 3 & 3 & 3 & 13 & 14 & 15 & 0.66 & 0.64 & 0.56 & 2.27 & 2.35 & 2.46  \\ \hline
  200 & 40 & ${\cal N}(0,3)$ & 4 & 3 & 3 & 33 & 36 & 35 & 0.46 & 0.42 & 0.50 & 2.30 & 2.44 & 2.39  \\ \hline
   200 & 60 & ${\cal N}(0,3)$ & 4 & 3 & 3 & 50 & 55 & 54 & 0.37 & 0.33 & 0.41 & 0.75 & 0.82 & 0.80  \\ \hline
      \end{tabular}
 \end{center}
\label{Tabl1} 
}
\end{table}

\subsection{Grouped variables}
For two error distributions, standard Normal ${\cal N}(0,1)$ and Cauchy ${\cal C}(0,1)$, we compare the results of the proposed  adaptive elastic-net method  with those of \cite{Ciuperca-17a} by adaptive group LASSO quantile method. The tuning parameters are for the our method: $\lambda_{1n}=c_2c(\sigma+c)g n^{\big((1-c)/2+1-(1-c)(1+\gamma)/2\big)/2}$, $\lambda_{2n}=c_3c(\sigma+c)n^{1/2 -c/2-1/n}$, with $c=\log(g)/\log(n)$  and $\gamma=\max \big(1.225,$ $2c/(1-c)+2/n\big)$, where $c_2$, $c_3$ are positive constants, varied on a value grid such that  $\lambda_{1n}$, $\lambda_{2n}$ make criterion (\ref{crit})  minimal. We took  $S_n=1$.  
We first consider that each group  contains $p=2$ explanatory variables, with $|{\cal A}^0|=4$, more precisely, $\eb^0_1=(0.5,1)$, $\eb^0_2=(1,1)$, $\eb^0_3=(-1,0)$, $\eb^0_4=(-1.5,1)$.  In Table \ref{Tabl_p2}, we give the same indicators as in the case of the ungrouped variables for two penalized estimation methods:   adaptive group LASSO quantile (\textit{aGq}) method and  adaptive elastic-net group quantile (\textit{aEGq}) method (\ref{aeQn}), for this last the algorithm presented in subsection  \ref{Algorithm} being used. On the explanatory variables, the both methods give similar results in the detection of the true significant groups of  and of the true non significant groups, but the \textit{aEGq} estimations are more accurate. To better exemplify this last finding, we give in Table  \ref{Tabl_p2_med} the medians of the parameter estimations. 
\begin{table} 
{\scriptsize
\caption{\footnotesize  Simulation results for models with grouped variables:  $p=2$. Estimation methods: adaptive group LASSO quantile \textit{(aGq)},  adaptive elastic-net group quantile (\textit{aEGq)}. }
\begin{center}
\begin{tabular}{|c|cc|c|c|c|c|c|c|c|c|c|c|c|c|}\hline
 & & &  \multicolumn{4}{c|}{true param $\neq 0$} & \multicolumn{4}{c|}{true param $= 0$} &  \multicolumn{2}{c|}{ } & \multicolumn{2}{c|}{ }\\
  \cline{4-15} 
 n &g & $\varepsilon$  &\multicolumn{2}{c|}{median} & \multicolumn{2}{c|}{mean} & \multicolumn{2}{c|}{median} & \multicolumn{2}{c|}{mean} & \multicolumn{2}{c|}{sd$(\beta^0-\hat \beta)$} & \multicolumn{2}{c|}{$mean|Y-\hat Y|$} \\
 & &  & aGq & aEGq& aGq & aEGq & aGq & aEGq & aGq & aEGq & aGq & aEGq & aGq & aEGq \\ \hline
  50 & 5 & ${\cal N}(0,1)$   & 4 & 4 & 3.9 & 3.99 & 0 & 1 & 0.47 & 0.85 & 1.01 & 0.99 & 1.03 & 1.14 \\
   &   & ${\cal C}(0,1)$   & 4 & 4 & 3.5 &  3.9 & 1 & 1 & 0.51 & 0.64 &  0.98 & 0.93 & 7.50& 7.59 \\
   \hline
    100 & 10 & ${\cal N}(0,1)$   & 4 & 3 & 3.6 & 2.88 & 5 & 6 & 4.9 & 6 & 0.67 & 0.62 & 1.35 & 1.82 \\
   &   & ${\cal C}(0,1)$   & 3 & 3 & 2.8 &  2.4 & 6 & 6 & 5.29 & 6 &  0.65 & 0.61 & 6.53 & 6.84 \\
   \hline
      \end{tabular}
 \end{center}
\label{Tabl_p2} 
}
\end{table}

\begin{table} 
{\scriptsize
\caption{\footnotesize   Median of $\widehat{\eb}_{{\cal A}^0}$, for models with grouped variables: $p=2$. Comparison adaptive elastic-net group quantile \textit{(aEGq)} estimations with adaptive group LASSO quantile \textit{(aGq)} estimations. }
\begin{center}
\begin{tabular}{|c|c|c|c|c|c|c|c|}\hline
 n &g & $\varepsilon$  &method &  $\eb^0_1=(0.5, 1)$ & $\eb^0_2=(1, 1)$ & $\eb^0_3=(-1, 0)$ & $\eb^0_4=(-1.5, 1)$ \\ \hline
    &  & ${\cal N}(0,1)$  & aGq & $ (0.34, 0.84)$ & $(0.67,  0.70)$ & $ (-0.37, 0)$ &  $(-0.82, 0.19)$\\
    &   &    & aEGq & $ (0.60, 0.74)$ & $(0.74,  0.74)$ & $ (-0.51, -0.20)$ &  $(-0.61, 0.16)$\\
 \cline{3-8}
  50  & 5 & ${\cal C}(0,1)$  & aGq & $ (0.22, 0.75)$ & $(0.45,  0.52)$ & $ (-0.05, 0)$ &  $(-0.41, 0)$\\
    &   &    & aEGq & $ (0.45, 0.55)$ & $(0.55,  0.55)$ & $ (-0.26, -0.06)$ &  $(-0.40, 0.13)$\\ \hline
     &   & ${\cal N}(0,1)$  & aGq & $ (0.22, 0.72)$ & $(0.55,  0.57)$ & $ (-0.09, 0)$ &  $(-0.48, 0)$\\
    &   &    & aEGq & $ (0.19, 0.24)$ & $(0.29,  0.30)$ & $ (-0.51, -0.20)$ &  $(-0.61, 0.16)$\\
  \cline{3-8}
  100  & 10 & ${\cal C}(0,1)$  & aGq & $ (0.07, 0.53)$ & $(0.33,  0.38)$ & $ (0, 0)$ &  $(-0.06, 0)$\\
    &   &    & aEGq & $ (0.13, 0.17)$ & $(0.20,  0.20)$ & $ (0, 0)$ &  $(-0.10, 0.03)$\\
  \hline
      \end{tabular}
 \end{center}
\label{Tabl_p2_med} 
}
\end{table}

To conclude our numerical study, we consider that each group of explanatory variables contains  $p=5$ variables. The first four groups are significants, with $\eb^0_1=(0.5,1,1.5,1,0.5)$, $\eb^0_2=(1,1,1,1,1)$, $\eb^0_3=(-1,0,1,2,1.5)$, $\eb^0_4=(-1.5,1,0.5,0.5,0.5)$, other groups are not significant. The design is generated as for the case $p=2$. The results are presented in Table \ref{Tabl2}. The finding made for groups of two variables remain true: the \textit{aEGq} method provides more accurate parameter estimations than those obtained by \textit{aGq}  method.  \\

The simulations were performed on a computer with CPU 1.90 GHz and 4 GB RAM. Execution time for 1000 Monte Carlo replications for $n=50$, $g=5$ and for a given value of $\lambda_{1n}$, $\lambda_{2n}$ in criterion (\ref{crit}) is 2.4 minutes, a single Monte Carlo replication taking 0.99 seconds. 
\begin{table} 
{\scriptsize
\caption{\footnotesize Simulation results for models with grouped variables: $p=5$. Comparison adaptive elastic-net group quantile \textit{(aEGq)} estimations with adaptive group LASSO quantile \textit{(aGq)} estimations. }
\begin{center}
\begin{tabular}{|c|cc|c|c|c|c|c|c|c|c|c|c|c|c|}\hline
 & & &  \multicolumn{4}{c|}{true param $\neq 0$} & \multicolumn{4}{c|}{true param $= 0$} &  \multicolumn{2}{c|}{ } & \multicolumn{2}{c|}{ }\\
  \cline{4-15} 
 n &g & $\varepsilon$  &\multicolumn{2}{c|}{median} & \multicolumn{2}{c|}{mean} & \multicolumn{2}{c|}{median} & \multicolumn{2}{c|}{mean} & \multicolumn{2}{c|}{sd$(\beta^0-\hat \beta)$} & \multicolumn{2}{c|}{$mean|Y-\hat Y|$} \\
 & &  & aGq & aEGq & aGq & aEGq & aGq & aEGq & aGq & aEGq & aGq & aEGq & aGq & aEGq \\ \hline
 50 & 5 & ${\cal N}(0,1)$   & 4 & 4 & 3.8 & 4 & 1 & 1 & 0.58 & 0.84 & 1.02 & 0.77 & 1.71 & 4.20 \\
    &  & ${\cal C}(0,1)$   & 4 & 4 & 3.6 & 3.9 & 1 & 1 & 0.74 & 0.55 & 1.04 & 0.81 & 6.98 & 7.25 \\
  \hline
 100 & 10 & ${\cal N}(0,1)$   & 4 & 4 & 3.85 & 4 & 6 & 6&  5.33 & 6 &  0.78 & 0.65 & 1.69 & 3.09 \\
     &  & ${\cal C}(0,1)$ &    4 & 4 & 3.6 & 3.9 & 6 & 6 & 5.6 & 5.8 & 0.77 & 0.63 & 7.37 & 8.42 \\
  \hline
   200 & 20 & ${\cal N}(0,1)$   & 4 & 4 & 3.8 & 4 & 15 & 16 & 15.2 & 16 & 0.57 & 0.47 & 1.85 & 3.66 \\
    &  & ${\cal C}(0,1)$    & 4 & 4 & 3.96 & 3.96 & 16 & 16 & 15.5 & 16 & 0.57 & 0.50 & 44.3 & 45.1 \\
  \hline
      \end{tabular}
 \end{center}
\label{Tabl2} 
}
\end{table}

 \section{Proofs}
  \label{Section_Proofs}
  In the proof of Theorem \ref{theorem1_SP} we will use the following theorem.
\begin{theorem}
\label{Geyer}
(\cite{Geyer-96}). If ${\cal F}_n$ and ${\cal F}$ are two random lower semicontinuous convex functions from $\R^p$  to $\R \cup \{ \infty \}$ such that $ {\cal F}_n\overset{\cal L} {\underset{n \rightarrow \infty}{\longrightarrow}} {\cal F}$ and ${\cal F}$ has a unique minimizer, then $\argmin_{\eu \in \R^p} {\cal F}_n (\eu)\overset{\cal L} {\underset{n \rightarrow \infty}{\longrightarrow}} \argmin_{\eu \in \R^p}{\cal F} (\eu).$
\end{theorem}
In the proofs, for all $x,y \in \R$,    the following  identity  on $\rho_\tau$ will be used:
\begin{equation}
\label{rho}
\rho_\tau(x-y)- \rho_\tau (x)=y(\e1_{x<0} - \tau )+ \int^{\tau}_0 (\e1_{x \leq v} - \eu_{x \leq 0})dv.
\end{equation}

\noindent{\bf Proof of Theorem \ref{th_vconv_SP}} 
By definition of adaptive elastic-net group quantile estimator, we have that 
\[
\widehat{{\eb}^g} =\argmin_{\eb^g \in \R^{r_n}} n^{-1} [E_n(\eb^g) -E_n(\eb^0)].
\]
For $\eu \in \R^{r_n}$, with $\| \eu \|_2=1$, $B>0$ a constant, in order to show the theorem, we study: $n^{-1} [E_n(\eb^0+B (g/n)^{1/2} \eu) -E_n(\eb^0)]$ which is, taking into account condition (\ref{C1bisb}) that $\lambda_{1n} n^{(c-1)/2- \alpha \gamma} {\underset{n \rightarrow \infty}{\longrightarrow}} 0$, using also the proof of Theorem 3 of  \cite{Ciuperca-17a}, strictly bigger that 
\begin{equation}
\label{RLL} 
\begin{split}
B^2 f(0) g n^{-1} \big( n^{-1}  \sum^n_{i=1} \eu^t \eeX_i \eeX_i^t \eu \big) (1+o_{\PP}(1)) - B O_{\PP}\pth{  g n^{-1}}\\
 +\lambda_{2n} n^{-1 }\sum^g_{j=1} \big[ \|\eb^0_j + g n^{-1} B \eu_j \|_2^2 - \| \eb^0_j\|_2^2 \big],
\end{split}
\end{equation}
for $B$ large enough. On the other hand, using  the triangular inequality, we have
\begin{align*}
\frac{\lambda_{2n}}{n }\sum^g_{j=1} \big[ \|\eb^0_j +\frac{g}{n} B \eu_j \|_2^2 - \| \eb^0_j\|_2^2 \big]& \leq \frac{\lambda_{2n}}{n }\sum^g_{j=1} \big[\frac{g}{n} B^2 \| \eu_j\|^2_2 +2 \| \eb^0_j \|_2 \sqrt{\frac{g}{n}} B \| \eu_j \|_2 \big] \\
&\simeq \lambda_{2n} \big[ \pth{\frac{g}{n}}^2 B^2+\pth{\frac{g}{n}}^{3/2} B \big] ,
\end{align*}
which is, using condition  (\ref{C1bis}), of order $o\pth{B^2 f(0) g n^{-1} \pth{n^{-1}  \sum^n_{i=1} \eu^t \eeX_i \eeX_i^t \eu } +B O_{\PP}\pth{ g n^{-1}} }$.  Then, taking into account relation (\ref{RLL}), together with Assumption (A6), we obtain that for all  $\epsilon >0$, there exists $B_\epsilon$  large enough such that for any $n$ large enough:
\[
\PP \bigg[ \inf_{\substack{\eu \in \R^{r_n} \\ \| \eu\|_2=1} } E_n \big(\ebo+B_\epsilon \sqrt{ \frac{g}{n}} \eu\big) >E_n(\ebo) \bigg]>1 - \epsilon
\]
and the theorem follows.
 \hspace*{\fill}$\blacksquare$ \\

\noindent {\bf Proof of Theorem \ref{Theorem 2SPL}}\\
The proof follows the same general lines as the proof of Theorem 4 in \cite{Ciuperca-17a} Consequently, some calculation details   are omitted.\\
\textit{(i)}  Let be the following two sets of parameters ${\cal V}_g(\ebo) \equiv \big\{ \eb^g \in \R^{r_n}; \| \eb^g-\ebo\|_2 \leq B (g/n)^{1/2}\big\}$ and ${\cal W}_n \equiv \big\{\eb^g \in {\cal V}_g(\ebo) ; \|\eb_{\overline {{\cal A}^0}}\|_2 >0\big\}$. By Theorem \ref{th_vconv_SP}, the estimator  $\widehat{{\eb}^g}$ belongs to ${\cal V}_g(\ebo)$ with a probability converging to 1 as $n \rightarrow \infty$ and for $B$ large enough. For proving the theorem, we will first show that 
\begin{equation}
\label{Pbg}
\lim_{n \rightarrow \infty} \PP  \big[ \widehat{{\eb}^g} \in{\cal W}_n  \big]=0.
\end{equation}
For this, we  consider  the parameter  vector $\eb^g=(\eb_{{\cal A}^0},  \eb_{\overline {{\cal A}^0}}) \in {\cal W}_n$. Let be also another parameter vector  $ {\eb}^{(1)}=({\eb}_{{\cal A}^0}^{(1)}, {\eb}_{\overline {{\cal A}^0}}^{(1)}) \in {\cal V}_g(\ebo)$,   such that ${\eb}^{(1)}_{{\cal A}^0}=\eb_{{\cal A}^0}$ and $ {\eb}_{\overline {{\cal A}^0}}^{(1)}=\textbf{0}_{r_n-p |{\cal A}^0|}$.\\
For proving relation (\ref{Pbg}), let us consider the following difference:
\begin{equation}
\begin{split}
D_n (\eb^g, \eb^{(1)}) & \equiv n^{-1} E_n(\eb^g) -  n^{-1}E_n({\eb}^{(1)})      \\
& =  \frac{1}{n} \sum^n_{i=1} \big[\rho_\tau(Y_i - \eeX^t_i \eb^g)-\rho_\tau(Y_i - \eeX^t_i {\eb}^{(1)})\big]      \\
&\qquad + \frac{\lambda_{1n}}{n}\sum^g_{j=|{\cal A}^0|+1}  \widetilde{\omega}_{n;j} \| \eb_j\|_2+  \frac{ \lambda_{2n}}{n}\sum^g_{j=|{\cal A}^0|+1} \| \eb_j\|_2^2.
\end{split}
\label{Dnbb}
\end{equation}
 Using identity (\ref{rho}), we have:
 \begin{align*}
 \sum^n_{i=1} \big[\rho_\tau(Y_i - \eeX^t_i \eb^g)-\rho_\tau(Y_i - \eeX^t_i {\eb}^{(1)})\big]&=  \sum^n_{i=1} (\eb^g-{\eb}^{(1)})^t \eeX_i \big[ \e1_{Y_i - \eeX_i^t {\eb}^{(1)} \leq 0}-\tau \big] \\
 &+\sum^n_{i=1} \int^{\eeX_i^t (\eb^g-{\eb}^{(1)})}_0 \big[\e1_{Y_i - \eeX_i^t {\eb}^{(1)} \leq v} -\e1_{Y_i - \eeX_i^t {\eb}^{(1)} \leq 0} \big] dv \\
 &  \equiv T_{1n}+T_{2n}.
\end{align*}
For $T_{1n}$, using assumption (A3), since $\eb^g \in {\cal V}_g(\ebo)$ and together with that $f$, $f'$ are bounded in a neighbourhood of 0 by assumption (A1),  we have, taking also into account assumption (A2), that,  $\E[T_{1n}]= \sum^n_{i=1} (\eb^g-\eb^{(1)})^t \eeX_i \big[F(\eeX_i^t ({\eb}^{(1)}-\ebo))-F(0) \big]=O\big(n \| \eb^g-\eb^{(1)} \|_2 \|\eb^{(1)}-\ebo\|_2  \big)=O\big(n \| \eb^g-\eb^{(1)} \|_2^2 \big)$. By similar arguments, using assumptions (A1), (A2), (A3), we have for the variance $\Var[T_{1n}]\leq \E[T^2_{1n}]= \sum^n_{i=1} \big((\eb^g-\eb^{(1)})^t \eeX_i \big)^2 =O\big(n \| \eb^g-\eb^{(1)}\|^2_2 \big)$. Then, by the Law of Large Numbers, we deduce,  
\begin{equation}
\label{ET1n}
T_{1n}= C n \| \eb^g-\eb^{(1)}\|^2_2 \big(1+o_{\PP}(1) \big)
\end{equation}
 For $T_{2n}=\sum^n_{i=1} \int^{\eeX_i^t (\eb^g-\eb^{(1)})}_0 \big[\e1_{\varepsilon_i  \leq v+\eeX_i^t (\eb^{(1)}-\eb^{0})} -\e1_{\varepsilon_i  \leq \eeX_i^t (\eb^{(1)}-\eb^{0})} \big] dv$,  since $\eb^g \in {\cal V}_g(\ebo)$,  we have, by a Taylor expansion, that, 
 \begin{equation*}
  \E[T_{2n}]=\sum^n_{i=1} \int^{\eeX_i^t (\eb^g-\eb^{(1)})}_0 \big[ F(v+\eeX_i^t (\eb^{(1)}-\eb^{0}))-F(\eeX_i^t (\eb^{(1)}-\eb^{0}) \big]dv .
 \end{equation*}
 Since the derivative $f'$ is bounded in a neighbourhood of 0, taking into account Assumption (A3), we have:
\begin{equation}
\label{ET2n}
  \E[T_{2n}]= O\big( \sum^n_{i=1} \big(\eeX_i^t (\eb^g-\eb^{(1)}) \big)^2\big) =O\big(n \| \eb^g-\eb^{(1)}\|_2^2 \big).
 \end{equation}
For the variance of $T_{2n}$, since the errors $\varepsilon_i$ are independent, we have:
  \begin{align*}
  \Var[ T_{2n}]&= \sum^n_{i=1} \Var \bigg[  \int^{\eeX_i^t (\eb^g-\eb^{(1)})}_0 \big[\e1_{\varepsilon_i  \leq v+\eeX_i^t (\eb^{(1)}-\eb^{0})} -\e1_{\varepsilon_i  \leq \eeX_i^t (\eb^{(1)}-\eb^{0})} \big]  dv \bigg]\\
  & = \sum^n_{i=1} \E \bigg[  \int^{\eeX_i^t (\eb^g-\eb^{(1)})}_0 \big( \big[\e1_{\varepsilon_i  \leq v+\eeX_i^t (\eb^{(1)}-\eb^{0})} -\e1_{\varepsilon_i  \leq \eeX_i^t (\eb^{(1)}-\eb^{0})} \big] \\
   & \qquad \qquad \qquad -\big[F(v+\eeX_i^t (\eb^{(1)}-\eb^{0})) -F(\eeX_i^t (\eb^{(1)}-\eb^{0})) \big] \big) dv \bigg]^2
  \\
 & \leq \sum^n_{i=1} \E \bigg[ \bigg| \int^{\eeX_i^t (\eb^g-\eb^{(1)})}_0 \big(  \big[\e1_{\varepsilon_i  \leq v+\eeX_i^t (\eb^{(1)}-\eb^{0})} -\e1_{\varepsilon_i  \leq \eeX_i^t (\eb^{(1)}-\eb^{0})} \big] -\big[F(v+\eeX_i^t (\eb^{(1)}-\eb^{0})) \\
  & \qquad \qquad \qquad -F(\eeX_i^t (\eb^{(1)}-\eb^{0})) \big] \big) dv \bigg| \bigg]  2 \big| \eeX_i^t (\eb^g-\eb^{(1)})\big|\\
&  \leq 2 \sum^n_{i=1} \int^{\eeX_i^t (\eb^g-\eb^{(1)})}_0  \big[F(v+\eeX_i^t (\eb^{(1)}-\eb^{0})) -F(\eeX_i^t (\eb^{(1)}-\eb^{0})) \big] dv\\
 & \qquad \qquad \qquad  \cdot 2 \max_{1 \leqslant l \leqslant n} \| \eeX_l\|_2 \| \eb^g-\eb^{(1)}\|_2.
  \end{align*}
  Taking into account assumptions (A1)-(A3) we obtain: $ \Var T_{2n}] =o(\E[T_{2n}] )$. Hence, taking also into account relation (\ref{ET2n}), by Bienaymé-Tchebychev inequality we obtain: $T_{2n}=C n \| \eb^g-\eb^{(1)}\|^2_2 \big(1+o_{\PP}(1) \big)$. Then, together with relation (\ref{ET1n}), we obtained that, $T_{1n}+T_{2n}=C n \| \eb^g-\eb^{(1)}\|^2_2 \big(1+o_{\PP}(1) \big)$. 
  Thus, returning to relation (\ref{Dnbb}), we get,
\begin{equation} 
 D_n (\eb^g, \eb^{(1)}) =  \displaystyle{C \| \eb^g - \eb^{(1)} \|_2^2(1+o_{\PP}(1)) + \frac{\lambda_{1n}}{n} \sum^g_{j=|{\cal A}^0|+1}  \pth{\frac{g}{n}}^{(1-\gamma)/2}+ \frac{ \lambda_{2n}}{n}\sum^g_{j=|{\cal A}^0|+1}  \frac{g}{n}}.
 \label{DD}
\end{equation}
From where,
\begin{eqnarray}
\frac{D_n(\eb^g, {\eb}^{(1)})}{\| \eb^g -{\eb}^{(1)} \|_2}& \geq &C \| \eb^g -{\eb}^{(1)} \|_2 (1+o_{\PP}(1))+  \frac{\lambda_{1n}}{n}     O_{\PP} \pth{ \pth{\frac{g}{n}}^{-\gamma/2}}+\frac{\lambda_{2n}}{n} \frac{g}{n} \nonumber \\
&=&O_{\PP}\bigg( \pth{\frac{g}{n}}^{1/2}+  \frac{\lambda_{1n}}{n}  \pth{\frac{g}{n}}^{-\gamma/2}+\frac{\lambda_{2n}}{n} \frac{g}{n} \bigg)\nonumber \\
&\geq& O_{\PP}\bigg(\frac{\lambda_{1n}}{n}  \pth{\frac{g}{n}}^{-\gamma/2} \bigg) =O_{\PP}\big(\lambda_{1n}n^{\gamma(1-c)/2-1} \big).
\label{Dn1}
\end{eqnarray}
We have similarly to relation (\ref{DD}) that $
D_n(\eb^0,\eb^{(1)})=C \|\big(\ebo - \eb^{(1)} \big)_{{\cal A}^0} \|_2^2 (1+o_{\PP}(1))$. 
From where, 
\begin{equation}
\label{Dn2}
\frac{D_n(\eb^0,\eb^{(1)})}{\| \ebo - \eb^{(1)} \|_2}=O_{\PP} \big(({g}/{n})^{1/2} \big)=O_{\PP}\big(n^{(c-1)/2} \big).
\end{equation}
Then, taking into account condition (\ref{C1bis2}), we obtain for relation  (\ref{Dn1},) by (\ref{Dn2}),
\[
\frac{D_n(\eb^g, {\eb}^{(1)})}{\| \eb^g -{\eb}^{(1)} \|_2}  \gg \frac{D_n(\eb^0,\eb^{(1)})}{\| \ebo - \eb^{(1)} \|_2},
\]
which implies relation (\ref{Pbg}). Relation (\ref{Pbg}) involves, with a probability converging to 1, as $n \rightarrow \infty$, that, for $\big\| \eb^g-\eb^0  \big\|_2 \leq C ( g/n )^{1/2}$, we have,
\begin{equation}
\label{eq15}
E_n(\eb^g_{{\cal A}^0}, \textbf{0}_{r_n-p|{\cal A}^0|})=\min_{\|\eb^g_{\overline {{\cal A}^0}}\|_2 \leq C \sqrt{\frac{g}{n}}} E_n(\eb^g_{{\cal A}^0}, \eb^g_{\overline {{\cal A}^0}}).
\end{equation}
On the other hand, by Assumption (A5) and Theorem \ref{th_vconv_SP} we have:   $\lim_{n \rightarrow \infty}\PP \big[   \min_{j \in {\cal A}^0} \| \widehat{{\eb}^g}_{j} \|_2 > 0  \big] = 1$. This relation together with relation (\ref{eq15}) imply claim \textit{(i)}.\\
\hh \textit{(ii)} By claim \textit{(i)}, we have with a probability converging to one as $n \rightarrow \infty$, that $\widehat{{\eb}^g}$ is under the form $\eb^0 +( {g}n^{-1})^{1/2} \ed $, with the $r_n$-dimensional vector  $\ed=(\ed_{{\cal A}^0},\textbf{0}_{r_n-p |{\cal A}^0|})$, $\ed_{{\cal A}^0}$ being a $(p |{\cal A}^0|)$-dimensional vector such that $\|\ed_{{\cal A}^0}\|_2 <C$. We consider the following difference:
\begin{equation}
\label{Enr}
\begin{array}{c}
\displaystyle{E_n \big(\ebo+ \sqrt{\frac{g}{n}} \ed\big)  - E_n(\ebo)=\frac{1}{n} \sum^n_{i=1} \big[\rho_\tau\big(Y_i-\eeX^t_i(\ebo+\sqrt{\frac{g}{n}} \ed)\big)-\rho_\tau(\varepsilon_i) } \big] \\
\displaystyle{+\frac{\lambda_{1n}}{n}\sum^{|{\cal A}^0|}_{j=1}  \widetilde{\omega}_{n;j} \big[ \| \eb^0_j+\sqrt{\frac{g}{n}} \ed_j\|_2 - \| \eb^0_j\|_2\big]+\frac{\lambda_{2n}}{n}\sum^{|{\cal A}^0|}_{j=1}   \big[ \| \eb^0_j+\sqrt{\frac{g}{n}} \ed_j\|_2^2 - \| \eb^0_j\|_2^2\big].}
\end{array}
\end{equation}
For relation (\ref{Enr}), we denote the following two sums ${\cal P}_1 \equiv n^{-1}\lambda_{1n} \sum^{|{\cal A}^0|}_{j=1}  \widetilde{\omega}_{n;j} \big[ \| \eb^0_j+(g/n)^{1/2} \ed_j\|_2 - \| \eb^0_j\|_2\big]$ and ${\cal P}_2 \equiv n^{-1}\lambda_{2n}\sum^{|{\cal A}^0|}_{j=1}   \big[ \| \eb^0_j+(g/n)^{1/2} \ed_j\|_2^2 - \| \eb^0_j\|_2^2\big]$.\\
Since $\| \widetilde {\eb^g} -\ebo \|_2=O_{\PP}((g/n)^{1/2}) $, we have ${\cal P}_1  = O_{\PP} \big(n^{-1}\lambda_{1n} |{\cal A}^0| (gn^{-1})^{1/2} \big)$. Then, using condition (\ref{C2bis1}), we obtain:
\begin{equation}
\label{P1gn}
{\cal P}_1/(g/n)=O_{\PP} \big( n^{-1}\lambda_{1n} |{\cal A}^0| (gn^{-1})^{-1/2} \big)= O_{\PP} \big( \lambda_{1n} |{\cal A}^0| n^{-(1+c)/2}  \big)=o_{\PP}(1).
\end{equation}
For the absolute value of ${\cal P}_2 $ we have by Cauchy-Schwarz inequality, that
\begin{align*}
\big| {\cal P}_2 \big| &\leq \frac{\lambda_{2n}}{n}\sum^{|{\cal A}^0|}_{j=1} \big(\frac{g}{n} \|\ed_j\|_2^2+2 \sqrt{\frac{g}{n}} \| \ed_j\|_2 \| \eb^0_j\|_2 \big)\\
 &  =O_{\PP} \big( \frac{\lambda_{2n}}{n}\sum^{|{\cal A}^0|}_{j=1} (g/n)^{1/2} \big)=  O_{\PP}(|{\cal A}^0| \lambda_{2n} n^{-1+(c-1)/2}).
\end{align*}
Then, using condition (\ref{C2bis2}), we have
\begin{equation}
\label{P2gn}
\frac{{\cal P}_2}{g/n} = O_{\PP}(\lambda_{2n} |{\cal A}^0| n^{-(1+c)/2})=o_{\PP}(1).
\end{equation}
For the first term of the right-hand side of (\ref{Enr}),  using the proof of Theorem 4 of  \cite{Ciuperca-17a}, we have: 
\begin{align}
\label{Pgn}
n^{-1} \sum^n_{i=1} \big[\rho_\tau\big(Y_i-\eeX^t_i(\ebo+\sqrt{\frac{g}{n}} \ed)\big)-\rho_\tau(\varepsilon_i) \big]&=\frac{1}{n} \sqrt{\frac{g}{n}} \sum^n_{i=1}\eeX^t_{i,{{\cal A}^0}} \ed_{{\cal A}^0} [\e1_{\varepsilon_i <0} - \tau]\nonumber \\
& \qquad +  \frac{f(0)}{2} \frac{g}{n} \ed^t_{{\cal A}^0} \eU_{n,{\cal A}^0} \ed_{{\cal A}^0}(1+o_{\PP}(1)).
\end{align}
Comparing (\ref{Pgn}) with (\ref{P1gn}) and (\ref{P2gn}) we deduce that minimizing (\ref{Enr}) amounts to minimize (\ref{Pgn}), which have that minimizer:
\begin{equation*}
\label{ad}
 \sqrt{\frac{g}{n}} \ed_{{\cal A}^0}=-\frac{1}{n} \frac{1}{f(0)} \eU^{-1}_{n,{\cal A}^0} \big( \sum^n_{i=1} \eeX_{i,{ {\cal A}^0}} (\e1_{\varepsilon_i \leq 0}-\tau)\big).
\end{equation*}
On the other hand, by Central Limit Theorem (CLT) for  independent random variable sequences , we have
\begin{equation*}
\label{cv}
\sqrt{n} f(0) \frac{\eu^t (\widehat{{\eb}^g}_{{\cal A}^0}- \eb^0_{{\cal A}^0})}{\sqrt{\tau(1- \tau) (\eu^t \eU^{-1}_{n,{\cal A}^0}  \eu)}}  \overset{\cal L} {\underset{n \rightarrow \infty}{\longrightarrow}} {\cal N}(0,1).
\end{equation*}
Since,   $\widehat{{\eb}^g}_{{\cal A}^0}- \eb^0_{{\cal A}^0}= (g/n)^{1/2} \ed_{{\cal A}^0}$, claim \textit{(ii)} follows.
   \hspace*{\fill}$\blacksquare$ \\

\noindent{\bf Proof of Remark \ref{theorem1_SP}}
Let us consider the vector  $\eu = (\eu_1, \cdots, \eu_g) \equiv \sqrt{n}(\eb^g- \ebo) \in \R^r $. Then $
Y_i-\eeX_i^t \eb^g= - n^{-1/2} \eeX_i^t \eu +\varepsilon_i$. We consider also the following random process:
\begin{align}
L_n(\eu)& \equiv \sum^n_{i=1} \big[ \rho_\tau \big( \varepsilon_i -n^{-1/2}  \eeX^t_i \eu \big)- \rho_\tau(\varepsilon_i) \big]+\lambda_{1n}\sum^g_{j=1}  \widetilde{\omega}_{n;j} \big[ \big\|\eb^0_j +n^{-1/2} \eu_j \big\|_2 - \| \eb^0_j \|_2  \big] \nonumber \\
 &  \hspace{5.7cm} + \lambda_{2n}\sum^g_{j=1} \big[ \big\|\eb^0_j +n^{-1/2} \eu_j  \big\|_2^2 - \| \eb^0_j \|_2^2  \big] \nonumber \\
  &=[\bfz^t_n \eu+B_n(\eu)]+\lambda_{1n}\sum^g_{j=1}  \widetilde{\omega}_{n;j} \big[ \big\|\eb^0_j +n^{-1/2} \eu_j \big\|_2 - \| \eb^0_j \|_2  \big] \frac{\sqrt{n}}{\sqrt{n}}  \nonumber \\
 &  \hspace{5.7cm}+ \lambda_{2n}\sum^g_{j=1} \big[ \big\|\eb^0_j +n^{-1/2} \eu_j  \big\|_2^2 - \| \eb^0_j \|_2^2  \big]  \nonumber\\
 &\equiv S_{1n} +S_{2n} +S_{3n},
 \label{S123n}
 \end{align}
with $\bfz_n  \equiv     n^{-1/2} \sum^n_{i=1} \eeX_i  \big( (1-\tau) \e1_{\varepsilon_i <0}- \tau \e1_{\varepsilon_i \geq 0}\big)$  and $B_n(\eu)  \equiv  \sum^n_{i=1} \int^{n^{-1/2}\eeX^t_i \eu }_0 [\e1_{\varepsilon_i < v} - \e1_{\varepsilon_i < 0} ]d v$. As in the proof of Theorem 1 of \cite{Ciuperca-17a} we have:  
\begin{equation}
\label{Wb}
S_{2n} \overset{\PP} {\underset{n \rightarrow \infty}{\longrightarrow}} \sum^g_{j=1} W(\eb^0_j, \eu), \quad  \textrm{with}  \quad W(\eb^0_j, \eu) \equiv \left\{
\begin{array}{lll}
0, & & \textrm{ if } \eb^0_j \neq \textbf{0}_p , \\
0, & &  \textrm{ if } \eb^0_j = \textbf{0}_p \textrm{ and } \eu_j=\textbf{0}_j, \\
\infty , & &  \textrm{ if } \eb^0_j = \textbf{0}_p \textrm{ and } \eu_j \neq \textbf{0}_j .
\end{array}
\right.
 \end{equation}
By Cauchy-Schwarz  inequality, we have: $\big|\|\ex - \ey \|^2_2 - \| \ex \|^2_2 \big|\leq \|\ey \|^2_2 +2 \|\ex\|_2 \| \ey \|_2$. Taking  $\ex=\eb^0_j$, $\ey = n^{-1/2} \eu_j$, using also condition (\ref{C1}), we obtain  for $S_{3n}$:
\begin{equation}
\label{S3n}
| S_{3n}| \leq \lambda_{2n} \sum^g_{j=1} \big[ n^{-1} \| \eu_j\|^2_2 +2 \| \eb^0_j\|_2 n^{-1/2} \|\eu_j\|_2 \big] = C \lambda_{2n} n^{-1/2} (1+o(1)) {\underset{n \rightarrow \infty}{\longrightarrow}} 0.
\end{equation}
For $S_{1n}$,  using Assumptions (A1), (A6),   we have,   by the CLT, that,   $ \bfz_n^t \eu \overset{\cal L} {\underset{n \rightarrow \infty}{\longrightarrow}} \bfz^t \eu$, with $\bfz$ a  $r$-dimensional random vector of law ${\cal N}(\textbf{0}_r, \tau(1-\tau)\eU )$. Taking into account assumption (A1) that the density $f$ has a bounded derivative in the neighborhood of 0, we have 
\begin{equation}
\label{EBn}
\E[B_n(\eu)] {\underset{n \rightarrow \infty}{\longrightarrow}} \frac{1}{2} f(0) \eu^t \eU \eu.
\end{equation}
On the other hand, 
\begin{align*}
\Var [B_n(\eu)] &=\sum^n_{i=1} \E \big[ \int^{n^{-1/2}\eeX^t_i \eu }_0 \big( \big( \e1_{\varepsilon_i < v} - \e1_{\varepsilon_i < 0} \big) -\big(F(v)-F(0)\big) \big) dv  \big]^2 \\
&\leq  \sum^n_{i=1} \E \big[\int^{n^{-1/2}\eeX^t_i \eu }_0  \big| \big( \e1_{\varepsilon_i < v} - \e1_{\varepsilon_i < 0} \big) -\big(F(v)-F(0)\big) \big| dv\big]  \cdot  2  n^{-1/2} \big|\eeX^t_i \eu  \big| \\
&\leq 2   n^{-1/2} \max_{1\leqslant i \leqslant n} \big|\eeX^t_i \eu  \big|  \sum^n_{i=1} \E \big[\int^{n^{-1/2}\eeX^t_i \eu }_0  \big( \big| \e1_{\varepsilon_i < v} - \e1_{\varepsilon_i < 0} \big| + \big|F(v)-F(0)\big| \big) dv\big]  \\
&= 4  n^{-1/2} \max_{1\leqslant i \leqslant n} \big|\eeX^t_i \eu  \big|  \sum^n_{i=1}  \int^{n^{-1/2}\eeX^t_i \eu }_0  \big(F(v)-F(0)\big)   dv\\
 &   =4 \E[B_n(\eu)]  n^{-1/2} \max_{1\leqslant i \leqslant n} \big|\eeX^t_i \eu  \big| {\underset{n \rightarrow \infty}{\longrightarrow}} 0,
 \end{align*}
by relation (\ref{EBn}) and Assumption (A6). Then, taking into account relation (\ref{EBn}), we obtain by Bienaym\'e-Tchebichev inequality that $B_n(\eu) \overset{\PP} {\underset{n \rightarrow  \infty}{\longrightarrow}} 2^{-1} f(0) \eu^t \eU \eu$.\\
Thus, taking into account relations (\ref{S123n}), (\ref{Wb}), (\ref{S3n}) and the fact that $ \bfz_n^t \eu \overset{\cal L} {\underset{n \rightarrow \infty}{\longrightarrow}} \bfz^t \eu$, together with $B_n(\eu) \overset{\PP} {\underset{n \rightarrow  \infty}{\longrightarrow}}  2^{-1} f(0) \eu^t \eU \eu$, we have:
\begin{equation}
\label{Lnu}
 L_n(\eu) \overset{\cal L} {\underset{n \rightarrow \infty}{\longrightarrow}} \bfz^t \eu +  \frac{1}{2} f(0) \eu^t \eU \eu +\sum^g_{j=1}W(\eb^0_j,\eu).
 \end{equation}
Since $L_n(\eu) $ and the right-hand side of relation (\ref{Lnu}) are convex in $\eu$, we have by Theorem \ref{Geyer} that: 
\begin{equation}
\label{argLn}
\argmin_{\eu \in \R^r} L_n(\eu) \overset{\cal L} {\underset{n \rightarrow \infty}{\longrightarrow}} \argmin_{\eu \in \R^r} \big(  \bfz^t \eu +  \frac{1}{2} f(0) \eu^t \eU \eu +\sum^g_{j=1}W(\eb^0_j,\eu) \big) .
\end{equation}
 On the other hand, we have $\widehat{{\eu}^n} \equiv \argmin_{\eu \in \R^r} L_n(\eu) \overset{denoted} {=}(\widehat{\eu}_{1n}, \widehat{\eu}_{2n})$, with $\widehat{\eu}_{1n}$  the first $(p|{\cal A}^0| )$ elements of $\widehat{{\eu}^n}$ and $\widehat{\eu}_{2n}$ the next  $(r-p|{\cal A}^0| )$ elements of $\widehat{{\eu}^n}$. Taking into account relations (\ref{argLn}), (\ref{Lnu}) and (\ref{Wb}) we obtain that $\widehat{\eu}_{2n} \overset{\PP} {\underset{n \rightarrow \infty}{\longrightarrow}} \textbf{0}_{r-p|{\cal A}^0|} $ and $\widehat{\eu}_{1n} \overset{\cal L} {\underset{n \rightarrow \infty}{\longrightarrow}} {\cal N}\big(\textbf{0}_{p|{\cal A}^0|}, \tau(1-\tau) f^{-2}(0)\eU^{-1}_{{\cal A}^0} \big)$. The proof of the theorem is thus complete. 
 \hspace*{\fill}$\blacksquare$


\end{document}